\documentclass[12pt]{amsart}
\input{amssymb.sty}
\usepackage[arrow, matrix]{xy}
\title{Weakly Exact von Neumann Algebras}
\author{Narutaka OZAWA}
\address{Department of Mathematical Sciences,
University of Tokyo, Komaba, 153-8914}
\email{narutaka@ms.u-tokyo.ac.jp}
\date{February 15, 2003}
\thanks{The author was supported by JSPS
Postdoctoral Fellowships for Research Abroad}

\newtheorem{thm}{Theorem}
\newtheorem{prop}[thm]{Proposition}
\newtheorem{cor}[thm]{Corollary}

\theoremstyle{definition}
\newtheorem{defn}[thm]{Definition}
\newcommand{\G}{\Gamma}
\newcommand{\e}{\varepsilon}
\newcommand{\p}{\varphi}
\newcommand{\B}{{\mathbb B}}

\newcommand{\M}{{\mathbb M}}
\newcommand{\N}{{\mathbb N}}  
 
\newcommand{\cb}{\mathrm{cb}} 
\newcommand{\hh}{{\mathcal H}} 
 
\newcommand{\id}{\mathrm{id}} 
\newcommand{\ip}[1]{\langle#1\rangle} 
\begin{document}
\begin{abstract}
The theory of exact $C^*$-algebras was introduced by 
Kirchberg and has been influential in recent development 
of $C^*$-algebras. 
A fundamental result on exact $C^*$-algebras is 
a local characterization of exactness. 
The notion of weakly exact von Neumann algebras was 
also introduced by Kirchberg. 
In this paper, we give a local characterization of 
weak exactness. 
As a corollary, we prove that a discrete group is 
exact if and only if its group von Neumann algebra is 
weakly exact. 
The proof naturally involves the operator space duality. 
\end{abstract}
\maketitle
\section{Introduction}
The theory of exact $C^*$-algebras was introduced and studied 
intensively by Kirchberg. 
We refer to \cite{wassermann} for general information on exact 
$C^*$-algebras. 
This influential theory has been playing a significant r{\^o}le 
in recent development of $C^*$-algebras. 
It is particularly important in the classification 
of $C^*$-algebras \cite{kicm}\cite{rordam} and 
in the theory of noncommutative topological entropy 
\cite{voiculescu}\cite{brown}\cite{stormer}. 
Hence it is natural to explore an analogue of this notion 
for von Neumann algebras. 

Throughout this paper, we mean by $J\triangleleft B$ 
a closed two sided ideal $J$ in a unital $C^*$-algebra $B$ 
and denote the quotient map by $Q\colon B\to B/J$. 
Although it is not essential, we assume for simplicity that 
all $C^*$-algebras, except ideals, are unital. 
We denote by $\otimes$ the algebraic tensor product 
and by $\otimes_{\min}$ the minimal (or spatial) tensor product. 
Recall that a $C^*$-algebra $A$ is said to be exact if 
\[
(A\otimes_{\min}B)/(A\otimes_{\min}J)=A\otimes_{\min}(B/J)
\]
for any $J\triangleleft B$, 
or equivalently, for any $*$-representation 
$\pi\colon A\otimes_{\min}B\to\B(\hh)$ 
with $A\otimes J\subset\ker\pi$, 
the induced $*$-representation 
$\tilde{\pi}\colon A\otimes(B/J)\to\B(\hh)$ 
is continuous w.r.t.\ the minimal tensor norm. 
Since von Neumann algebras are not exact unless 
they are subhomogeneous, the correct notion of exactness 
for von Neumann algebras shall be 
the weak exactness introduced in \cite{kicm}. 
Thus, for a von Neumann algebra $A=M$, 
we modify the above definition and require 
the $*$-representation $\pi$ to be left normal, i.e., 
the restriction of $\pi$ to $M$ is normal. 
\begin{defn}\label{defwex}
A von Neumann algebra $M$ is said to be weakly exact if 
for any $J\triangleleft B$ and 
any left normal $*$-representation 
\[
\pi\colon M\otimes_{\min}B\to\B(\hh)
\]
with $M\otimes J\subset\ker\pi$, 
the induced $*$-representation 
\[
\tilde{\pi}\colon M\otimes(B/J)\to\B(\hh)
\]
is continuous w.r.t.\ the minimal tensor norm. 
\end{defn}
Kirchberg \cite{kicm} proved that 
a von Neumann algebra $M$ is weakly exact if 
it contains a weakly dense exact $C^*$-algebra $A$.
Indeed, 
by Theorem 4 (10) in \cite{kirchberg}, 
the exact $C^*$-algebra $A$ is locally reflexive 
\cite{ab}\cite{eh} i.e., the left normal embedding 
$A^{**}\otimes_{\min} C\subset(A\otimes_{\min}C)^{**}$ 
is continuous for any $C$. 
Let $J\triangleleft B$, $\pi$ and $\tilde{\pi}$ 
be given as in Definition~\ref{defwex}.
Since $A$ is exact, $\tilde{\pi}$ is continuous 
on $A\otimes_{\min}(B/J)$. 
By the local reflexivity of $A$, 
$\tilde{\pi}$ extends to a left normal $*$-representation 
$\rho$ on $A^{**}\otimes_{\min}(B/J)$. 
But, since $\tilde{\pi}$ is left normal on $M\otimes(B/J)$, 
it coincides with $\rho$ on $M\otimes(B/J)$ 
and hence is continuous on $M\otimes_{\min}(B/J)$. 
We will prove a partial converse of this result and 
a local characterization of weak exactness. 
Before stating the theorem, we recall that 
an operator system $S$ is exact if 
\[
(S\otimes_{\min} B)/(S\otimes_{\min} J)=S\otimes_{\min}(B/J)
\] 
isometrically for any $J\triangleleft B$. 
We note that exact operator systems are characterized locally 
\cite{kirchberg}\cite{pisierex} and they are locally reflexive 
(Corollary 4.8 in \cite{eor}). 
The term ``u.c.p.'' stands for ``unital completely positive''. 

\begin{thm}\label{thm}
Let $M$ be a von Neumann algebra with a separable predual. 
Then, the following are equivalent. 
\begin{enumerate}
\item\label{wex1}
The von Neumann algebra $M$ is weakly exact.
\item\label{wex2}
For any finite dimensional operator system $E$ in $M$, 
there exist nets of u.c.p.\ maps 
$\p_i\colon E\to\M_{n(i)}$ 
and $\psi_i\colon\p_i(E)\to M$ such that 
the net $\{\psi_i\circ\p_i\}$ converges to $\id_E$ 
in the point-ultrastrong topology. 
\item\label{wex3}
There exist an exact operator system $S$ and 
normal u.c.p.\ maps 
$\p\colon M\to S^{**}$ and $\psi\colon S^{**}\to M$ 
such that $\psi\circ\p=\id_M$. 
\end{enumerate}
\end{thm}

We note that the implications 
$\ref{wex1}\Leftrightarrow\ref{wex2}\Leftarrow\ref{wex3}$ 
hold without separability assumption on $M$. 
It would be interesting to know whether 
we can take $S$ in the condition~\ref{wex3} 
to be an exact $C^*$-algebra 
(and moreover $\psi$ to be a normal $*$-homomorphism). 
Kirchberg \cite{kicm} asked whether every von Neumann algebra is 
weakly exact. 
The following corollary solves this problem negatively. 
(See \cite{gromov} for existence of non-exact groups.)
It is left open whether the ultrapower $R^\omega$ of 
the hyperfinite type $\mathrm{II}_1$ factor $R$ is weakly exact or not. 

\begin{cor}\label{cor}
Let $\G$ be a discrete group. 
Then, the group von Neumann algebra $L\G$ is weakly exact 
if and only if $\G$ is exact. 
\end{cor}

We recall that a discrete group $\G$ is exact if and only if 
its reduced group $C^*$-algebra $C^*_{\mathrm{r}}\G$ 
is exact \cite{kw}. 
Since $C^*_{\mathrm{r}}\G$ is weakly dense in $L\G$, 
the ``if'' part of the corollary follows from 
Kirchberg's theorem mentioned above. 
It is not hard to show that weak exactness passes to 
a von Neumann subalgebra which is the range of 
a normal conditional expectation. 
Hence, every von Neumann subalgebra of a weakly exact 
\emph{finite} von Neumann algebra is again weakly exact. 
It follows that if $\G$ is exact and $\Lambda$ is not, 
then $L\Lambda\not\hookrightarrow L\G$. 
Moreover, exactness is a measure equivalence invariant. 
\section{Proofs}
It is well-known that exactness for a $C^*$-algebra 
is equivalent to the property $C'$ (cf.\ \cite{ab}\cite{eh}). 
We prove a von Neumann algebra analogue of this fact. 
Recall that for von Neumann algebras $M$ and $N$, 
a $*$-representation $\pi$ of $M\otimes N$ is said to be 
\emph{binormal} if both $\pi|_M$ and $\pi|_N$ are normal. 

\begin{prop}\label{prop}
A von Neumann algebra $M$ is weakly exact if and only if 
for any $C^*$-algebra $B$ and any left normal $*$-representation 
$\pi\colon M\otimes_{\min}B\to\B(\hh)$, 
the binormal extension 
$\hat{\pi}\colon M\otimes B^{**}\to\B(\hh)$ 
is continuous w.r.t.\ the minimal tensor norm. 
\end{prop}
\begin{proof}
We first prove the ``if'' part. 
Let $J\triangleleft B$ and $\pi\colon M\otimes_{\min}B\to\B(\hh)$ 
be as in Definition~\ref{defwex}. 
Let $p$ be the central projection which supports 
the normal $*$-homomorphism 
$Q^{**}\colon B^{**}\to(B/J)^{**}$ so that 
we may identify $(B/J)^{**}$ with $pB^{**}$. 
We denote the canonical inclusion by 
\[
\psi\colon B/J\to (B/J)^{**}=pB^{**}\subset B^{**}.
\]
By the assumption, $\pi$ extends to a binormal 
$*$-homomorphism $\hat{\pi}$ on $M\otimes_{\min}B^{**}$. 
Since $M\otimes J\subset\pi$, 
we have $\hat{\pi}(1\otimes p)=1$. 
We claim that $\tilde{\pi}$ coincides with 
$\hat{\pi}\circ(\id\otimes\psi)$ and hence 
continuous on $M\otimes_{\min}(B/J)$. 
Indeed, for any $a\in M$ and $x\in B$, we have 
\[
(\hat{\pi}\circ(\id\otimes\psi))(a\otimes Q(x))
=\hat{\pi}(a\otimes px)
=\pi(a\otimes x)=\tilde{\pi}(a\otimes Q(x)).
\]

We next prove the ``only if'' part. 
Let a $C^*$-algebra $B$ and a left normal $*$-representation 
$\pi\colon M\otimes B\to\B(\hh)$ be given.   
For a directed set $I$, we set
\[
B_I=\{(x(i))_{i\in I} \in \prod_{i\in I}B :
\mbox{ultrastrong$^*$-}\lim_{i\in I} x(i)
\mbox{ exists in }B^{**}\}.
\]
Since adjoint operation and multiplication are jointly
ultrastrong$^*$-continuous on a bounded subset,
$B_I$ is a C$^*$-algebra and the map
\[
\sigma\colon B_I\ni(x(i))_{i\in I}
\mapsto\mbox{ultrastrong$^*$-}\lim_{i\in I} x(i)\in B^{**}
\]
is a continuous $*$-homomorphism.
Because of Kaplansky's density theorem,
we may assume that $\sigma$ is surjective.
Consider the $*$-representation 
$\rho\colon M\otimes B_I\to\B(\hh)$ 
defined by 
\[
\rho(\sum_{k=1}^n a_k\otimes (x_k(i))_{i\in I})=
\mbox{ultrastrong$^*$-}\lim_i\pi(\sum_{k=1}^n a_k\otimes x_k(i)).
\]
We observe that $\rho$ is continuous on 
$M\otimes_{\min}B_I\subset\prod_{i\in I}(M\otimes_{\min}B)$.
Since $M$ is weakly exact, $\rho$ is left normal and 
$M\otimes\ker\sigma\subset\ker\rho$, the induced $*$-representation 
$\tilde{\rho}$ is continuous on $M\otimes_{\min}B^{**}$.
It is not too hard to see that $\hat{\pi}=\tilde{\rho}$ and hence 
$\hat{\pi}$ is continuous on $M\otimes_{\min}B^{**}$. 
\end{proof}

The advantage of the above formulation 
is that $B$ in the statement need not be a $C^*$-algebra.
We make this point more precise. 
Let $X$ be a $C^*$-algebra or an operator space. 
Let $p\in M^{**}$ be the central support of 
the identity representation of $M$ so that 
we may identify $M$ with $pM^{**}$. 
Then, the canonical binormal embedding 
$M^{**}\otimes X^{**}\subset(M\otimes_{\min}X)^{**}$ 
gives rise to a (non-unital) binormal embedding 
\[
\Theta_X\colon M\otimes X^{**} = pM^{**}\otimes X^{**} 
\hookrightarrow (M\otimes_{\min}X)^{**}, 
\]
i.e., for $z=\sum a_k\otimes x_k\in M\otimes X$, we have 
$\Theta_X(z)=\sum pa_k\otimes x_k$. 
Now, assume that $M$ is weakly exact and 
let $B$ be a C$^*$-algebra. 
Since $\Theta_B$ is a binormal $*$-homomorphism which 
is continuous on $M\otimes_{\min}B$, it follows from 
Proposition~\ref{prop} that $\Theta_B$ is continuous 
(isometric) on $M\otimes_{\min}B^{**}$. 
This implies that $\Theta_X$ is isometric 
on $M\otimes_{\min}X^{**}$ for any operator space $X$. 
Indeed, when $X\subset B$, we have the commuting diagram; 
\[\xymatrix{
 M\otimes X^{**} \ar[r]^{\Theta_X}\ar@{^{(}->}[d]
 & (M\otimes_{\min} X)^{**}\ar@{^{(}->}[d]\\
 M\otimes_{\min} B^{**} \ar[r]^{\Theta_B}
 & (M\otimes_{\min} B)^{**}
}\]
where the bottom and the vertical inclusions 
are all isometric. 
It follows that for any $z\in M\otimes X^{**}$ 
with $\|z\|_{\min}\le 1$, there exists a net 
$\{z_i\}$ in $M\otimes X$ with $\|z_i\|_{\min}\le1$
which converges to $\Theta_X(z)\in (M\otimes X)^{**}$ 
in the weak$^*$-topology. 
We observe that the net $\{z_i\}$ converges to $z$ 
in the $\sigma(M\otimes X^{**}, M_*\otimes X^*)$-topology. 

We use the operator space duality to restate the above result. 
Recall that the operator space duality says that 
for $z=\sum a_k\otimes x_k\in M\otimes X^{**}$ and 
the corresponding (finite rank) map $T_z\colon X^*\to M$ 
defined by 
\[
T_z\colon X^*\ni f\mapsto\sum \ip{f,x_k}\,a_k\in M,
\]
we have $\|T_z\|_{\cb}=\|z\|_{\min}$
(cf.\ \cite{erbook}\cite{pisierbook}). 
We note that $T_z$ is weak$^*$-continuous
if and only if $z$ comes from $M\otimes X$ 
(rather than $M\otimes X^{**}$).
Thus we arrive at the following corollary. 
The term ``c.c.'' stands for ``completely contractive''. 
\begin{cor}\label{wcp}
Let $M$ be a weakly exact von Neumann algebra and 
$X$ be an operator space. 
Then, for any finite rank c.c.\ map $\p\colon X^*\to M$, 
there exists a net $\{\p_i\}$ of weak$^*$-continuous 
finite rank c.c.\ maps $\p_i\colon X^*\to M$ 
which converges to $\p$ in the point-ultraweak topology. 
\end{cor}
We are now in position to prove the theorem.
\begin{proof}[Proof of Theorem \ref{thm}]
\textbf{Ad}\ref{wex1}$\Rightarrow$\ref{wex2}: 
Let $M\subset\B(\hh)$ be a weakly exact von Neumann algebra 
and $E\subset M$ be a finite dimensional operator system. 
Fix an increasing sequence $\{\hh_k\}$ of finite dimensional 
subspaces of $\hh$ with a dense union and 
denote by $\Phi_k\colon E\to\B(\hh_k)=\M_{n(k)}$ 
the compression corresponding to $\hh_k$. 
We note that $\Phi_k$ is a linear isomorphism onto 
$E_k=\Phi_k(E)$ if $k$ is sufficiently large. 
Consider a complete isometry given by
\[
\Phi\colon E\ni a\mapsto(\Phi_k(a))_k\in\prod_k E_k.
\] 
We note that $\prod E_k\subset\prod\M_{n(k)}$ 
is an ultraweakly closed operator subspace (and hence 
is a dual operator space). 
There exists a right inverse of $\Phi$ given by 
\[
\Psi\colon\prod_k E_k\ni(x_k)_k\mapsto
\mbox{ultraweak-}\lim_{k\to\omega}\Phi_k^{-1}(x_k)
\in E\subset M, 
\]
where $\omega$ is a fixed free ultrafilter on $\N$. 
Indeed, $\Psi$ is a well-defined c.c.\ map 
with $\Psi\circ\Phi=\id_E$. 
By Corollary~\ref{wcp}, there exists a net of 
ultraweakly continuous c.c.\ maps $\psi_i\colon\prod E_k\to M$ 
such that $\lim_i\psi_i=\Psi$ 
in the point-ultraweak topology. 
Since each $\psi_i$ is ultraweakly continuous, 
if we set 
\[
\p_m\colon E\ni a\mapsto (\Phi_1(a),\ldots,\Phi_m(a),0,0,\ldots)
\in\prod_k E_k,
\]
then we have 
\[
\lim_i\lim_m\psi_i\circ\p_m
=\lim_i\psi_i\circ\Phi=\Psi\circ\Phi=\id_E.
\] 
We note that 
$\p_m(E)\subset\bigoplus_{k=1}^m E_k
\subset\bigoplus_{k=1}^m\M_{n(k)}$.
Passing to a convex combination, we may assume that 
the convergence is w.r.t.\ the point-ultrastrong topology.
This proves the condition~\ref{wex2}
except that $\psi_i$'s may not be u.c.p. 
One can fix this problem as follows. 
First, using a standard cut and paste method, one may assume that 
$\psi_i$ is approximately unital in norm. 
Then, one invokes Theorem 2.5 in \cite{eh}. 

\textbf{Ad}\ref{wex2}$\Rightarrow$\ref{wex3}: 
Let $M$ be a weakly exact von Neumann algebra with a separable predual 
and $A$ be a ultraweakly-dense norm-separable C$^*$-subalgebra in $M$. 
We denote by $|\,\cdot\,|$ a seminorm which defines the ultrastrong
topology on the unit ball of $M$. 
Using the condition~\ref{wex2} recursively, we can construct 
sequences of finite dimensional operator systems 
and connecting u.c.p.\ maps 
\[\xymatrix{
E_1 \ar[rr]^{\theta_1} \ar[dr]_{\p_1} & & 
E_2 \ar[rr]^{\theta_2} \ar[dr]_{\p_2} & & 
E_3 \ar[rr]^{\theta_3} \ar[dr]_{\p_3} & & \cdots\\
& F_1 \ar[rr]^{\sigma_1} \ar[ur]^{\psi_1} &  & 
F_2 \ar[rr]^{\sigma_2} \ar[ur]^{\psi_2} & & 
F_3 \ar[rr]^{\sigma_3} \ar[ur]^{\psi_3} & & \cdots
}\]
such that 
\begin{enumerate}
\item 
$E_1\subset E_2\subset\cdots\subset M$ 
and the norm-closure of $\bigcup E_k$ contains $A$, 
\item
$F_n\subset\M_{n(k)}$ for some $n(k)\in\N$,
\item
the diagram commutes and 
$|\theta_k(a)-a|<2^{-k}\|a\|$ for $a\in E_k$.
\end{enumerate}
We define the operator system $S$ 
as the inductive limit of $(F_k,\sigma_k)$; 
let 
\[
S=Q(\tilde{S})\subset\prod\M_{n(k)}/\bigoplus\M_{n(k)}
\]
where $Q\colon\prod\M_{n(k)}\to\prod\M_{n(k)}/\bigoplus\M_{n(k)}$ 
is the quotient map and $\tilde{S}$ is the norm closure of 
\[
\{ (x_k)_k\in\prod_{k=1}^\infty\M_{n(k)} : 
x_k\in F_k\mbox{ and }x_{k+1}=\sigma_k(x_k)
\mbox{ eventually}\}\subset\prod_{k=1}^\infty\M_{n(k)}.
\]
We note that $S$ coincides with the quotient of $\tilde{S}$ 
by the complete $M$-ideal $\bigoplus\M_{n(k)}$. 
Since $\tilde{S}$ is exact and locally reflexive, so is $S$. 
Let $\Phi\colon\bigcup E_k\to S^{**}$ be a cluster point 
of the sequence 
\[
\Phi_m\colon\bigcup E_k\ni a\mapsto 
Q((\sigma_{k-1}\circ\cdots\circ\sigma_m\circ\p_m
(a))_{k=m+1}^\infty)\in S
\] 
and let $\tilde{\psi}\colon\tilde{S}\to M$ be a cluster point 
of the sequence 
\[
\tilde{\psi}_m\colon\tilde{S}\ni (x_k)_k\mapsto\psi_m(x_m)\in M.
\]
In both cases, cluster points are taken w.r.t.\ 
the point-weak$^*$ topology.
It is easy to see that $\tilde{\psi}=\psi\circ Q$ for some 
u.c.p.\ map $\psi\colon S\to M$. 
The unique weak$^*$-continuous extension of $\psi$ on $S^{**}$ 
is still denoted by $\psi$. 
Then, for any $a\in\bigcup E_n$, we have that 
\begin{align*}
|a-\psi\circ\Phi(a)| 
&\le \limsup_m |a-\psi\circ\Phi_m(a)|\\
&= \limsup_m |a-\tilde{\psi}((\sigma_{k-1}\circ\cdots
 \circ\sigma_m\circ\p_m(a))_{k=m+1}^\infty)|\\
&\le \limsup_m\limsup_k |a-\psi_k\circ\sigma_{k-1}\circ
 \cdots\circ\sigma_m\circ\p_m(a)|\\
&= \limsup_m\limsup_k |a-\theta_k\circ\cdots\circ\theta_m(a)|\\
&\le \limsup_m\limsup_k \sum_{j=m}^k 2^{-j}\|a\|=0.
\end{align*}
It follows that $\psi\circ\Phi$ is the identity map on $\bigcup E_k$. 
We first extend $\Phi$ on the norm closure of $\bigcup E_k$ 
by norm continuity and then further extend $\Phi|_A$ 
to $\bar{\Phi}\colon A^{**}\to S^{**}$ by ultraweak continuity. 
Then, $\psi\circ\bar{\Phi}\colon A^{**}\to M$ is a normal 
u.c.p.\ map such that $\psi\circ\bar{\Phi}|_A=\id_A$. 
Hence, the restriction $\p$ of $\bar{\Phi}$ to $M\subset A^{**}$ 
satisfies $\psi\circ\p=\id_M$. 

\textbf{Ad}\ref{wex3}$\Rightarrow$\ref{wex1}: 
We omit the proof because it is almost same as that of Kirchberg's 
theorem cited in the remarks preceding Theorem~\ref{thm}. 
\end{proof}
\begin{proof}[Proof of Corollary \ref{cor}]
Let $\G$ be a discrete group such that $L\G$ is weakly exact.
To prove exactness of $\G$, we use the criterion 
given in \cite{exact}. 
Let a finite subset $\mathcal{E}\subset\G$ and $\e>0$ be given. 
By Arveson's extension theorem, Theorem~\ref{thm} and its proof 
(or Paulsen's trick), there exist a finite subset 
$\mathcal{F}\subset\G$ and a u.c.p.\ map 
$\psi\colon\B(\ell_2\mathcal{F})\to\B(\ell_2\G)$ 
such that if we denote by 
$\p\colon L\G\to\B(\ell_2\mathcal{F})$ 
the compression, then $\psi\circ\p(\lambda(s))\in L\G$ and 
$\|\lambda(s)-\psi\circ\p(\lambda(s))\|_2<\e$
for every $s\in\mathcal{E}$. 
As in \cite{gk}, we consider the positive definite kernel 
defined by 
\[
u(s,t)=\ip{\psi\circ\p(\lambda(st^{-1}))\delta_t,\delta_s}
\]
for $s,t\in\G$. 
We have that $u(s,t)\neq0$ only if 
$st^{-1}\in\mathcal{F}\mathcal{F}^{-1}$ 
and that $\|u(s,t)-1\|<\e$ if $st^{-1}\in\mathcal{E}$. 
This proves exactness of $\G$.
\end{proof}


\begin{thebibliography}{EOR}
%
\bibitem[AB]{ab}
R. J. Archbold and C. J. K. Batty. 
\textit{$C^*$-tensor norms and slice maps.} 
J. London Math. Soc. (2) \textbf{22} (1980), no. 1, 127--138. 
%
\bibitem[Br]{brown}
N. P. Brown, 
\textit{Topological entropy in exact $C^*$-algebras.} 
Math. Ann. \textbf{314} (1999), no. 2, 347--367. 
%
\bibitem[EH]{eh}
E. G. Effros and U. Haagerup. 
\textit{Lifting problems and local reflexivity for $C^*$-algebras.} 
Duke Math. J. \textbf{52} (1985), no. 1, 103--128.
%
\bibitem[EOR]{eor}
E. G. Effros, N. Ozawa and Z.-J. Ruan, 
\textit{On injectivity and nuclearity for operator spaces.} 
Duke Math. J. \textbf{110} (2001), 489--521.
%
\bibitem[ER]{erbook}
E. G. Effros and Z.-J. Ruan. 
\textit{Operator spaces.} 
London Mathematical Society Monographs. 
New Series, \textbf{23}. 
The Clarendon Press, Oxford University Press, New York, 2000. 
%
\bibitem[Gr]{gromov}
M. Gromov. 
\textit{Random walk in random groups.} 
Preprint IHES (2001).
%
\bibitem[GK]{gk}
E. Guentner and J. Kaminker. 
\textit{Exactness and the Novikov conjecture.} 
Topology \textbf{41} (2002), no. 2, 411--418. 
%
\bibitem[Ki1]{kicm}
E. Kirchberg. 
\textit{Exact $C^*$-algebras, tensor products, 
and the classification of purely infinite algebras.} 
Proceedings of the International Congress of Mathematicians, 
Vol. 1, 2 (Z{\"u}rich, 1994), 943--954, Birkh{\"a}user, Basel, 1995.
%
\bibitem[Ki2]{kirchberg}
E. Kirchberg. 
\textit{On subalgebras of the CAR-algebra.} 
J. Funct. Anal. \textbf{129} (1995), no. 1, 35--63.
%
\bibitem[KW]{kw}
E. Kirchberg and S. Wassermann. 
\textit{Permanence properties of $C^*$-exact groups.} 
Doc. Math. \textbf{4} (1999), 513--558. 
%
\bibitem[Oz]{exact}
N. Ozawa.
\textit{Amenable actions and exactness for discrete groups.} 
C. R. Acad. Sci. Paris Sér. I Math. \textbf{330} (2000), 
no. 8, 691--695.
%
\bibitem[Pi1]{pisierex}
G. Pisier. 
\textit{Exact operator spaces.} 
Recent advances in operator algebras (Orl{\'e}ans, 1992). 
Ast{\'e}risque No. \textbf{232} (1995), 159--186. 
%
\bibitem[Pi2]{pisierbook}
G. Pisier. 
\textit{Introduction to Operator Space Theory}. 
Cambridge University Press (2003). 
%
\bibitem[R{\o}]{rordam}
M. R{\o}rdam. 
\textit{Classification of nuclear, simple $C^*$-algebras.} 
``Classification of nuclear $C\sp *$-algebras. 
Entropy in operator algebras'', 1--145, 
Encyclopaedia Math. Sci., \textbf{126}, Springer, Berlin, 2002. 
%
\bibitem[St]{stormer} 
E. St{\o}rmer, 
\textit{A survey of noncommutative dynamical entropy.} 
``Classification of nuclear $C\sp *$-algebras. 
Entropy in operator algebras'', 147--198, 
Encyclopaedia Math. Sci., \textbf{126}, Springer, Berlin, 2002.
%
\bibitem[Vo]{voiculescu}
D. Voiculescu. 
\textit{Dynamical approximation entropies 
and topological entropy in operator algebras.} 
Comm. Math. Phys. \textbf{170} (1995), no. 2, 249--281. 
%
\bibitem[Wa]{wassermann}
S. Wassermann, 
\textit{Exact $C^*$-algebras and related topics.} 
Lecture Notes Series \textbf{19}. 
Seoul National University, Research Institute of Mathematics, 
Global Analysis Research Center, Seoul, 1994.
%
\end{thebibliography}
\end{document}